\input amstex
\input amsppt.sty
\NoBlackBoxes
\nologo
\topmatter
\title A correction to the paper \linebreak
``Log Abundance Theorem for Threefolds"
\endtitle
\author Kenji Matsuki
\endauthor
\NoPageNumbers
\leftheadtext{}
\rightheadtext{}
\endtopmatter
\document

Dr. Qihong Xie of Tokyo University points out that in Chapter 6 of the paper

``Log Abundance Theorem for Threefolds" by Sean Keel, Kenji Matsuki, and James
McKernan, Duke Math. J. Vol. 75 No. 1 (1994), 99-119,
there are several crucial mistakes and misleading statements.  

The original manuscripts
of the second author, on which the entire Chapter 6 was based upon, did not contain
these mistakes, and the statements were different from the ones in the paper.  The
calamity somehow crept into the paper in the process of collaboration in order to adapt,
simplify and modify the arguments in the manuscripts.

\vskip.1in

To be more precise and concrete, we list the mistakes and problematic statements that
Xie [22] communicated to us below:

\vskip.1in

Theorem 6.1: In ``Moreover" statement in case (2) of the assertion, the crucial
assumption of $D_1 \cdot D_2 \cdot ... \cdot D_{n-1}$ being not numerically trivial is
missing.  Furthermore, the asserted inequality is
different from the one in the original manuscripts mentioned above. 

\vskip.1in

The absence of the assumption also affects the statements of Corollaries 6.2 and 6.3,
both of which did not exist in the original manuscripts.

\vskip.1in

Corollary 6.2: The paper intended to apply Theorem 6.1 with $D_{n-1} = D_n = -
K_X$ and the fact that there is no family of rational curves $C$ covering $X$ with $-
K_X
\cdot C = 0$.  The correct statement with the extra assumption will be found at the end
of this correction.

\vskip.1in

Corollary 6.3: The statement is correct, but the proof needs a little more care taking
the extra assumption into consideration.
\vskip.1in

Lemma 6.4: The second part of the assertion ``and there is a fixed positive number
$\epsilon$ such that $f \leq g$ implies $f(a) \leq g(a)$, for all $0 \leq a \leq
\epsilon$ and any $f$ and $g$ in $P$." is a sheer nonsense.  For example, if we take $P
= \{f_n, g_n;n \in {\Bbb N}, f_n =  - n^2a + 1, g_n = - (n+1)^2 a + 2\}$, then
$P$ is a subset of ${\Bbb Q}[a]$ such that the degree, the coefficients, and the
denominator of any coefficient are bounded from above.  Obviously we have $f_n \leq
g_n$  according to the definition of the paper.  However, as $f_n - g_n = (2n + 1)a -
1$, there is no fixed positive number
$\epsilon$ such that $f_n(a) \leq g_n(a)$ for all $0 \leq a \leq \epsilon$ and for all
$n \in {\Bbb N}$.

\vskip.1in

Lemma 6.5: In ``Furthermore" part of the assertion, the crucial
assumption of $D_1 \cdot D_2 \cdot ... \cdot D_{n-1}$ being not numerically trivial is
missing.  Moreover, the asserted inequality is different from the one
in the original manuscripts mentioned above.

\vskip.1in

Proof of Lemma 6.5: For the assertion that the Harder-Narasimhan
filtration of $E$ with respect to $H_1 \cdot H_2 \cdot ... \cdot H_{n-1}$ is
independent of $a$ when $a$ is small enough (and positive), the paper simply
attributes the proof to Lemma 6.4, which is found to be a nonsense as above.  In the
original manuscripts, there was no Lemma 6.4 and the argument there for the stability
(independence) of the Harder-Narasimhan filtration is different and more explicit than
the one in the paper, as will be presented in this correction.

For the proof of the inequality, the paper says, ``(Now, if we set \linebreak
$D =
\frac{1}{k}c_1(F) -
\frac{1}{(r - k)}c_1(E/F)$, then $D \cdot D_1 \cdot D_2 \cdot ... D_{n-1} = 0$,) and
so, by applying the Hodge index theorem on $S$ and taking the limit, $D^2 \cdot D_1
\cdot D_2
\cdot ... \cdot D_{n-2} \leq 0$.  Using this, it is an easy calculation to deduce the
Bogomolov inequality for $E$."  Not only that the author (Kenji Matsuki) of this
correction does not understand this part of the argument, but also, since it seems that
it does not use the crucial assumption of $D_1 \cdot D_2 \cdot ... \cdot D_{n-1}$ being
not numerically trivial in any way, it raises a serious question to the validity of the
proof (and also to the asserted inequality itself).  In the original manuscripts,
the inequality asserted is different from the one in the paper, and the explicit
computation to derive the inequality shows clearly one place where the assumption of
$D_1 \cdot D_2 \cdot ...
\cdot D_{n-1}$ being not numerically trivial is used in an essential way.

\vskip.1in

In this correction, therefore, we rewrite Chapter 6 presenting the arguments closer to
and more faithful to the original manuscripts, in order to clarify these mistakes and
the statements.  With the corrected arguments and statements, we do estabilsih the
required inequality involving the second Chern class to complete the arguments for the
log abundance conjecture (for threefolds) in the case $\nu = 2$.  The rest of the proof
remains valid without any change.

Although the original paper is a collaboration of the above three authors, the
second author, Kenji Matsuki, is solely responsible for this correction, as well as for
the mistakes and misleading statements in the published paper which called for the
necessity of this correction.  

\vskip.1in

We thank Dr. Qihong Xie for bringing the attention of the author to the needs of this
correction, and Profs. Keiji Oguiso and Yoichi Miyaoka for valuable help.

\vskip.2in

\proclaim{New Chapter 6. Interaction of Bogomolov stability and Mori's bend and break
method}\endproclaim

\proclaim{6.1. Brief review of the outline}\endproclaim  We briefly recall the essence
of the log abundance conjecture for threefolds.

\vskip.1in

Let $(X,\Delta)$ be a pair consisting of a normal projective (or more generally
complete) variety $X$ of dimension 3 (over an algebraically closed field of
characteristic zero) and an effective ${\Bbb Q}$-divisor $\Delta$, called a boundary
divisor, such that it has only kawamata log terminal (or more generally log canonical)
singularities and that $K_X + \Delta$ is nef.  The conjecture asserts that the linear
system
$|m(K_X + \Delta)|$ is base point free for some $m \in {\Bbb N}$.

The basic and rough outline of our proof for the conjecture goes as follows:

Firstly, we show that $|m(K_X + \Delta)|$ is non-empty for some $m \in {\Bbb N}$.  

This
is done in Theorem 1.2 of Keel-Matsuki-McKernan[8].  The argument follows closely the
one by Miyaoka[13,14,15] for the existence of a member of $|mK_X|$, for some $m \in
{\Bbb N}$, where $X$ is a minimal threefold, as the first step to prove the usual
abundance conjecture for threefolds (cf.Kawamata[5]).  In our log(arithmic) version, a
special attention has to be paid to the case where $X$ has a structure of a uniruled
variety with
$\kappa(X,K_X) = -
\infty$ and where the addition of an effective ${\Bbb Q}$-divisor should push up
$\kappa(X,K_X + \Delta)$ to $0$ or more.

Secondly, we analyze the cases according to the numerical Kodaira dimension $\nu(X,K_X
+ \Delta) = 0,1,2,3$.

Suppose $\nu = 0$.  Then the existence of a member implies that $m(K_X + \Delta) \sim
0$ for some $m \in {\Bbb N}$, verifying the log abundance.

Suppose $\nu = 3 = \dim X$.  Then the base point freeness theorem of
Kawamata-Reid-Shokurov (cf. Kawamata-Matsuda-Matsuki[7]) tells us that $|m(K_X +
\Delta)|$ is base point free for sufficiently large and divisible $m \in {\Bbb N}$,
verifying the log abundance.

In order to deal with the remaining cases where $\nu = 1, 2$, we follow the strategy of
Kawamata[5] applied to settle the usual abundance conjecture for threefolds in the cases
where
$\nu = 1,2$.

One of the main ingredients of his strategy is the application of the log MMP to find a
``better" (minimal) model than the original one, whose abundance statement, however,
will imply that of the original.  This part can be carried out similarly for our
log(arithmic) version.  (We use the same notation $(X,\Delta)$ for the ``better" model
by abuse of notation.)  See Chapter 5 of Keel-Matsuki-McKernan[8] for detail and the
precise meaning of the word ``better".

Once this is done, the proof for the case where $\nu = 1$ goes almost verbatim as the
one in Chapter 13 of Koll\'ar et al[9] for the usual abundance.

Another of the main ingredients of his strategy, to settle the case where $\nu = 2$,
is the non-negativity of the intersection $c_2 \cdot L$ where $L = mK_X$. 
This guarantees that there are enough sections to imply the abundance, via the
computation of the Riemann-Roch formula for $L$ and via the log abundance for
surfaces in an inductional scheme of the proof.  The non-negativity is a consequence of
the Bogomolov-Miyaoka inequality, which is implied by the generic semi-positivity of
$\Omega_X^1$ (which is equivalent to the generic semi-negativity of the tangent sheaf
$T_X$).  

In our log(arithmic) version, however, we do not have the generic semi-positivity of
$\Omega_X^1$, as $X$ may well be uniruled as an ambient space.  

This is the place where
we have to observe the interaction of the Bogomolov stability and Mori's bend
and break method more closely, as we do here in Chapter 6.

Namely, if the tangent sheaf is not semi-negative, then there is some destabilizing
``positive" subsheaf.  Then by the method of Mori's bend and break, there exists a
family of rational curves $C$ along this subsheaf with the property $L \cdot C = 0$
where $L = m(K_X + \Delta)$.  Essentially the family of rational curves induces a fiber
space structure
$X
\rightarrow Y$, where the line bundle $L$ is the pull-back of another on $Y$.  Thus the
log abundance of $L$ on $X$ is reduced to that of this line bundle on $Y$, which has
dimension one less, and we are done.

If the tangent sheaf does not have such a destabilizing subsheaf, i.e., if it is
semistable, then we can derive the same non-negativity of the intersection $c_2 \cdot L$
via the modification of the Bogomolov-Miyaoka inequality, obtaining the same estimate
for the sections of $L$ as before.  This finishes the proof for the
case where $\nu = 2$.

\vskip.1in

In the following, we present this interaction of the Bogomolov stability and
Mori's bend and break method in detail.

\vskip.1in

\proclaim{6.2. Setting}\endproclaim Let $X$ be a normal projective variety of dimension
$d$ which has quotient singularities in codimension 2, i.e., there is a Zariski closed
subset $F$ of
$X$ with
$\roman{codim}_XF \geq 3$ such that $X \setminus F$ has only quotient singularities. 
(This is the case if, for example, $X$ has only (kawamata) log terminal
singularities.)

We remark then that, though the notion of ${\Bbb Q}$-varieties and ${\Bbb Q}$-sheaves
in the sense of Mumford[18] (cf. Matsusaka[11]Kawamata[5]Chapter 10 in Koll\'ar et
al[9]) is well-defined only over $X \setminus F$, for a ${\Bbb Q}$-sheaf on $X \setminus
F$, the Chern classes
$\hat{c}_1({\Cal E})$ and
$\hat{c}_2({\Cal E})$ are well-defined objects in $A^1(X) \otimes {\Bbb Q}$ and
$A^2(X) \otimes {\Bbb Q}$, respectively, via the isomorphisms
$$\align
A^1(X \setminus F) \otimes {\Bbb Q} &\overset{\sim}\to{\rightarrow} A^1(X)
\otimes {\Bbb Q}
\\
A^2(X \setminus F) \otimes {\Bbb Q} &\overset{\sim}\to{\rightarrow} A^2(X)
\otimes {\Bbb Q}.
\\
\endalign$$

(In Chapter 10 of Koll\'ar et al[9], the Chern classes for ${\Bbb Q}$-varieties are
explicitly defined by Megyesi only in dimension 2, where the condition in Mumford
[18] of the ``global cover" $\tilde{X}$ being Cohen-Macaulay is automatically
satisfied.  Being a non-expert in the field, we can only offer an explicit
reference without the dimension restriction of the Chern classes (the resolution
property) to the recent paper by Totaro[20].  The reader is also invited to look at
Vistoli[21], Eddidin-Hassett-Kresch-Vistoli[2] among others.)

We also remark that, for Cartier divisors $D_1, D_2, ... , D_{d-1}$ on $X$, the
intersection numbers
$\hat{c}_1({\Cal E}) \cdot D_1 \cdot D_2 \cdot ... \cdot D_{d-1}$,
$\hat{c}_1({\Cal E})^2 \cdot D_1 \cdot D_2 \cdot ... \cdot D_{d-2}$, and \linebreak
$\hat{c}_2({\Cal E}) \cdot D_1 \cdot D_2 \cdot ... \cdot D_{d-2}$ are well-defined. 
(Note that, by taking appropriate hyperplane sections, these intersection numbers can
be computed on normal surfaces with only quotient singularities.)

We write $K_X = \hat{c}_1(\Omega_X^1) = - \hat{c}_1(T_X) \in A^1(X)
\otimes {\Bbb Q}$, where $\Omega_X^1$ is the sheaf of 1-forms and $T_X$ is the tangent
sheaf of $X$, respectively.

For a torsion free coherent sheaf ${\Cal E}$ on $X$ and nef Cartier divisors $D_1, D_2,
... , D_{d-1}$, it is straightforward to define the (averaged) slope
$$\mu_{(D_1 \cdot D_2 \cdot ... \cdot
D_{d-1})}({\Cal E}) = \frac{\hat{c}_1({\Cal E}) \cdot D_1 \cdot D_2 \cdot ... \cdot
D_{d-1}}{\roman{rank}\ {\Cal E}}$$
and semistability of ${\Cal E}$ with respect to $(D_1 \cdot D_2 \cdot ... \cdot
D_{d-1})$, as well as to show the existence and the uniqueness of the Harder-Narasimhan
filtration.

\vskip.1in

\proclaim{6.3. Theorem} Let $D_1, D_2, ... , D_d$ be nef Cartier divisors
on $X$ such that \linebreak
$D_1 \cdot D_2 \cdot ... \cdot D_d = 0$ and $- K_X \cdot D_1
\cdot D_2
\cdot ... \cdot D_{d-1} \geq 0$.  

Then either

(1) $X$ is covered by a family of rational curves (of bounded degree) such that $D_d
\cdot C = 0$, or 

(2) the tangent sheaf $T_X$ is $(D_1 \cdot D_2 \cdot ... \cdot D_{d-1})$-semistable
and \linebreak
$\hat{c}_1(T_X) \cdot D_1 \cdot D_2 \cdot ... \cdot D_{d-1} = - K_X \cdot D_1 \cdot D_2
\cdot ... \cdot D_{d-1} = 0$.

Moreover, in case (2), if $D_1 \cdot D_2 \cdot ... \cdot D_{d-1}$ is not numerically
trivial, then we have the inequality
$$\hat{c}_2(T_X) \cdot D_1 \cdot D_2 \cdot ... \cdot D_{d-2} \geq
\frac{1}{2}\hat{c}_1(T_X)^2 \cdot D_1 \cdot D_2 \cdot ... \cdot D_{d-2}.$$ 
(In particular, if we assume further that $\hat{c}_1(T_X)^2 \cdot D_1 \cdot D_2 \cdot
...
\cdot D_{d-2} = K_X^2
\cdot D_1 \cdot D_2 \cdot ... \cdot D_{d-2} \geq 0$, then $\hat{c}_2(T_X) \cdot D_1
\cdot D_2 \cdot ... \cdot D_{d-2} \geq 0$.)

\endproclaim

The rest of the chapter is devoted to the proof of Theorem 6.3.

\vskip.1in

\proclaim{6.4. Proposition} Let ${\Cal E}$ be a torsion free coherent sheaf, $H$
an ample divisor, and $D_1, D_2, ... , D_{d-1}$ nef Cartier divisors on $X$.  Then the
Harder-Narasimhan filtration of ${\Cal E}$ with respect to $(H + nD_1) \cdot (H +
nD_2) \cdot ... \cdot (H + nD_{d-1})$ 
$$0 = {\Cal E}_{0,n} \subset {\Cal E}_{1,n} \subset ... \subset {\Cal E}_{k_n,n} =
{\Cal E}$$
stabilizes for $n$ sufficiently large.  That is to say, there exists a filtration
$$0 = {\Cal E}_0 \subset {\Cal E}_1 \subset ... \subset {\Cal E}_k =
{\Cal E}$$
and $n_o \in {\Bbb N}$ such that for $n \geq n_o$ we have
$$\CD
0 @. = @. {\Cal E}_0 @. \subset @. {\Cal E}_1 @. \subset @.\hskip.1in 
@....@.\hskip.1in  @.
\subset @. {\Cal E}_k = {\Cal E} \\
@| @.@| @.@| @.@.@.@.@.@|   \\
0 @. = @. {\Cal E}_{0,n} @. \subset @. {\Cal E}_{1,n} @. \subset @. @.... @.@. \subset
@. {\Cal E}_{k_n,n} = {\Cal E} \\
\endCD$$
\endproclaim

\demo{Proof}\enddemo 

Step 0. Let 
$$n^{d-1}A_{d-1} + n^{d-2}A_{d-2} + ... + nA_1 + A_0 = (H + nD_1) \cdot (H +
nD_2) \cdot ...
\cdot (H + nD_{d-1})$$
be the expansion according to the powers of $n$.  Remark the two facts.  The first,
which can be easily verified since $X$ has quotient singularities in codimension 2, is
that there exists
$q
\in {\Bbb N}$ such that for any subsheaf
${\Cal F}
\subset {\Cal E}$ we have
$$\mu_{A_{d-1}}({\Cal F}), \mu_{A_{d-2}}({\Cal F}), ... , \mu_{A_1}({\Cal F}),
\mu_{A_0}({\Cal F}) \in
\frac{1}{q
\cdot r!}{\Bbb Z}$$ where $r = \roman{rank}\ {\Cal E}$.  The second is that there
exists $b
\in {\Bbb N}$ such that for any subsheaf ${\Cal F} \subset {\Cal
E}$
$$\mu_{A_{d-1}}({\Cal F}), \mu_{A_{d-2}}({\Cal F}), ... , \mu_{A_1}({\Cal F}),
\mu_{A_0}({\Cal F}) \leq b.$$
In order to see the second assertion, one can for example argue as follows.  Take the
dual ${\Cal E}^{\vee} = Hom_{{\Cal O}_{X}}({\Cal E},{\Cal O}_X)$.  Choose a
sufficiently ample divisor $L$ such that ${\Cal E}^{\vee} \otimes L$ is globally
generated, i.e., we have a surjection ${\Cal O}_X^{\oplus m} \rightarrow {\Cal
E}^{\vee} \otimes L \rightarrow 0$.  Taking the dual sequence, we have an injection $0
\rightarrow {\Cal E}^{\vee\vee} \rightarrow L^{\oplus m}$.  Note that for any subsheaf
${\Cal F} \subset {\Cal E}$ there is an open subset $U \subset X$ with
$\roman{codim}_X(X \setminus U) \geq 2$ such that ${\Cal F}|_U$ and ${\Cal E}|_U$
are locally free and that ${\Cal E}|_U = {\Cal E}^{\vee\vee}|_U$.  The above injection
then induces another injection $0 \rightarrow {\Cal F}|_U
\rightarrow L^{\oplus m}|_U$ and hence one more injection $0 \rightarrow
\wedge^s{\Cal F}|_U \rightarrow \wedge^s(L^{\oplus m})|_U$ where $s = \roman{rank}\
{\Cal F}$.  Since $\roman{codim}_X(X \setminus U) \geq 2$, this implies that
$$\mu_{A_i}({\Cal F}) = \frac{1}{s}\hat{c}_1({\Cal F}) \cdot A_i \leq L \cdot A_i \leq
\max_i\{L \cdot A_i\} = b,$$
where $b = \max_i\{L \cdot A_i\}$ is independent of ${\Cal F}$.

\vskip.1in

Step 1. Let $m_{d-1} = \max\{\mu_{A_{d-1}}({\Cal F});{\Cal F} \subset {\Cal E}\}$.  The
existence of such $m_{d-1}$ is guaranteed by Step 0.  We claim that there is a maximal
subsheaf ${\Cal E}_{A_{d-1}}$ with $\mu_{A_{d-1}}({\Cal E}_{A_{d-1}}) = m_{d-1}$, i.e.,
any subsheaf ${\Cal G} \subset {\Cal E}$ with $\mu_{A_{d-1}}({\Cal G}) = m_{d-1}$ is
contained in ${\Cal E}_{A_{d-1}}$.  

In fact, let ${\Cal G}, {\Cal H} \subset {\Cal E}$
be two subsheaves with $\mu_{A_{d-1}}({\Cal G}) = \mu_{A_{d-1}}({\Cal H}) = m_{d-1}$. 
Consider the exact sequence
$$0 \rightarrow {\Cal G} \cap {\Cal H} \rightarrow {\Cal G} \oplus {\Cal H} \rightarrow
{\Cal G} + {\Cal H} \rightarrow 0$$
where ${\Cal G} \cap {\Cal H}, {\Cal G} + {\Cal H} \subset {\Cal E}$ are subsheaves and
where $(g,h) \in {\Cal G} \oplus {\Cal H} \mapsto g - h \in {\Cal G} + {\Cal H}$. 
We have
$$\align
&\frac{\roman{rank}\ ({\Cal G} \cap {\Cal H})\mu_{A_{d-1}}({\Cal G} \cap {\Cal H}) +
\roman{rank}\ ({\Cal G} + {\Cal H})\mu_{A_{d-1}}({\Cal G} +
{\Cal H})}{\roman{rank}\ ({\Cal G} \cap {\Cal H}) + \roman{rank}\ ({\Cal G} + {\Cal
H})} \\
&=
\mu_{A_{d-1}}({\Cal G}
\oplus {\Cal H}) \\
&= \frac{\roman{rank}\ ({\Cal G})\mu_{A_{d-1}}({\Cal G}) +
\roman{rank}\ ({\Cal H})\mu_{A_{d-1}}({\Cal H})}{\roman{rank}\ ({\Cal G})
+ \roman{rank}\ ({\Cal H})}\\
&= m_{d-1}\\
\endalign$$ 
Since $\mu_{A_{d-1}}({\Cal G} \cap {\Cal H}),\mu_{A_{d-1}}({\Cal G} +
{\Cal H}) \leq m_{d-1}$ by the maximality of $m_{d-1}$, we conclude
$$\mu_{A_{d-1}}({\Cal G} \cap {\Cal H}) = \mu_{A_{d-1}}({\Cal G} +
{\Cal H}) = m_{d-1}.$$
Now the existence of the maximal subsheaf ${\Cal E}_{A_{d-1}}$ follows from the
noetherian property of the coherent sheaf ${\Cal E}$.

\vskip.1in

Inductively we construct the subsheaves ${\Cal E}_{(A_{d-1},A_{d-2}, ... ,
A_{d-i})}$ for $i = 2, ... , d$.

\vskip.1in

Step $i$. Let $m_{d-i} = \max\{\mu_{A_{d-i}}({\Cal F});{\Cal F} \subset {\Cal E},
\mu_{A_{d-k}}({\Cal F}) = m_{d-k}, k = 1, ... , i-1\}$.  The existence of such $m_{d-i}$
is guaranteed by Step 0.  We claim that there is a maximal subsheaf ${\Cal
E}_{(A_{d-1},A_{d-2}, ... , A_{d-i})}$ with $\mu_{A_{d-k}}({\Cal E}_{(A_{d-1},A_{d-2},
... , A_{d-i})}) = m_{d-k}$ for $k = 1, ... , i$, i.e., any subsheaf ${\Cal G} \subset
{\Cal E}$ with $\mu_{A_{d-k}}({\Cal G}) = m_{d-k}$ for $k = 1, ... , i$ is contained in
${\Cal E}_{(A_{d-1},A_{d-2}, ... , A_{d-i})}$.

In fact, let ${\Cal G}, {\Cal H} \subset {\Cal E}$
be two subsheaves with $\mu_{A_{d-k}}({\Cal G}) = \mu_{A_{d-k}}({\Cal H}) = m_{d-k}$
for $k = 1, ... , i$.  Consider the exact sequence
$$0 \rightarrow {\Cal G} \cap {\Cal H} \rightarrow {\Cal G} \oplus {\Cal H} \rightarrow
{\Cal G} + {\Cal H} \rightarrow 0$$
where ${\Cal G} \cap {\Cal H}, {\Cal G} + {\Cal H} \subset {\Cal E}$ are subsheaves and
where $(g,h) \in {\Cal G} \oplus {\Cal H} \mapsto g - h \in {\Cal G} + {\Cal H}$.
Then in Step ($i-1$) we showed as an inductional hypothesis that
$$\mu_{A_{d-k}}({\Cal G} \cap {\Cal H}) = \mu_{A_{d-k}}({\Cal G} +
{\Cal H}) = m_{d-k} \text{\ for\ }k = 1, ... , i-1.$$
We have
$$\align
&\frac{\roman{rank}\ ({\Cal G} \cap {\Cal H})\mu_{A_{d-i}}({\Cal G} \cap {\Cal H}) +
\roman{rank}\ ({\Cal G} + {\Cal H})\mu_{A_{d-i}}({\Cal G} +
{\Cal H})}{\roman{rank}\ ({\Cal G} \cap {\Cal H}) + \roman{rank}\ ({\Cal G} + {\Cal
H})} \\
&=
\mu_{A_{d-i}}({\Cal G}
\oplus {\Cal H}) \\
&= \frac{\roman{rank}\ ({\Cal G})\mu_{A_{d-i}}({\Cal G}) +
\roman{rank}\ ({\Cal H})\mu_{A_{d-i}}({\Cal H})}{\roman{rank}\ ({\Cal G})
+ \roman{rank}\ ({\Cal H})}\\
&= m_{d-i}.\\
\endalign$$ 
Since $\mu_{A_{d-i}}({\Cal G} \cap {\Cal H}),\mu_{A_{d-i}}({\Cal G} +
{\Cal H}) \leq m_{d-i}$ by the maximality of $m_{d-i}$, we conclude
$$\mu_{A_{d-i}}({\Cal G} \cap {\Cal H}) = \mu_{A_{d-i}}({\Cal G} +
{\Cal H}) = m_{d-i}.$$
Now the existence of the maximal subsheaf ${\Cal E}_{(A_{d-1},A_{d-2}, ... ,A_{d-i})}$
follows from the noetherian property of the coherent sheaf ${\Cal E}$.  Note that
$${\Cal E}_{(A_{d-1},A_{d-2}, ... , A_{d-i})} \subset {\Cal E}_{(A_{d-1},A_{d-2}, ...
,A_{d-(i-1)})} \subset {\Cal E}.$$

\vskip.1in

By induction on the rank of ${\Cal E}$, it is enough to show that the first piece of
the Harder-Narasimhan filtration ${\Cal E}_{1,n}$ is stable for $n \gg 0$.

\vskip.1in

Final Step.  We claim that ${\Cal E}_{(A_{d-1},A_{d-2}, ... ,A_1,A_0)}$, which we
denote by ${\Cal E}_1$ for short, is the first piece ${\Cal E}_{1,n}$ of the
Harder-Narasimhan filtration with respect to \linebreak
$(H + nD_1)
\cdot (H + nD_2) \cdot ... \cdot (H + nD_{d-1})$ for $n \gg 0$.

\vskip.1in

Choose $n_1 \in {\Bbb N}$ such that for all $n \geq n_1$ we have
$$m_{d-1}n^{d-1} + m_{d-2}n^{d-2} + ... + m_1n + m_0 > \{m_{d-1} - \frac{1}{q \cdot
r!}\}n^{d-1} + bn^{d-2} + ... + bn + b.$$
For a subsheaf ${\Cal F} \subset {\Cal E}$, if $\mu_{A_{d-1}}({\Cal F}) \neq
m_{d-1}$, then we have for $n \geq n_1$
$$\align
&\mu_{(H + nD_1) \cdot (H + nD_2) \cdot ... \cdot (H + nD_{d-1})}({\Cal E}_1) \\
&= m_{d-1}n^{d-1} + m_{d-2}n^{d-2} + ... + m_1n + m_0 \\
&> \{m_{d-1} - \frac{1}{q \cdot
r!}\}n^{d-1} + bn^{d-2} + ... + bn + b \\
&\geq \mu_{A_{d-1}}({\Cal F})n^{d-1} + \mu_{A_{d-2}}({\Cal F})n^{d-2} + ... +
\mu_{A_1}({\Cal F})n + \mu_{A_0}({\Cal F})\\
&= \mu_{(H + nD_1) \cdot (H + nD_2) \cdot ... \cdot (H + nD_{d-1})}({\Cal F}). \\
\endalign$$
Therefore, in order for ${\Cal F}$ to have the maximal (averaged) slope with respect to
\linebreak
$(H + nD_1) \cdot (H + nD_2) \cdot ... \cdot (H + nD_{d-1})$ it is necessary
that
$\mu_{A_{d-1}}({\Cal F}) = m_{d-1}$.

\vskip.1in
 
Inductively, for $i = 2, ... , d$, choose $n_1 \leq n_2 \leq ... \leq n_i$ such
that for all $n \geq n_i$ we have
$$\align
&m_{d-1}n^{d-1} + m_{d-2}n^{d-2} + ... + m_{d-(i-1)}n^{d - (i-1)}\\
&\hskip2in + m_{d-i}n^{d-i} +
m_{d-(i+1)}n^{d-(i+1)} + ... + m_1n + m_0
\\ &> m_{d-1}n^{d-1} + m_{d-2}n^{d-2} + ... + m_{d-(i-1)}n^{d - (i-1)}\\
&\hskip2in + \{m_{d-i}
- \frac{1}{q \cdot r!}\}n^{d-i} + bn^{d-(i+1)} + ... + bn + b.\\
\endalign$$ 
For a subsheaf ${\Cal F}
\subset {\Cal E}$ with 
$$\mu_{A_{d-1}}({\Cal F}) = m_{d-1}, \mu_{A_{d-2}}({\Cal F}) = m_{d-2}, ... ,
\mu_{A_{d-(i-1)}}({\Cal F}) = m_{d-(i-1)},$$  if
$\mu_{A_{d-i}}({\Cal F}) \neq m_{d-i}$, then we have for $n \geq n_i$
$$\align
&\mu_{(H + nD_1) \cdot (H + nD_2) \cdot ... \cdot (H + nD_{d-1})}({\Cal E}_1) \\
&= m_{d-1}n^{d-1} + m_{d-2}n^{d-2} + ... + m_{d-(i-1)}n^{d - (i-1)}\\
&\hskip2in + m_{d-i}n^{d-i} +
m_{d-(i+1)}n^{d-(i+1)} + ... + m_1n + m_0
\\ 
&> m_{d-1}n^{d-1} + m_{d-2}n^{d-2} + ... + m_{d-(i-1)}n^{d - (i-1)}\\
&\hskip2in + \{m_{d-i}
- \frac{1}{q \cdot r!}\}n^{d-i} + bn^{d-(i+1)} + ... + bn + b.\\
&\geq \mu_{A_{d-1}}n^{d-1} + \mu_{A_{d-2}}n^{d-2} + ... + \mu_{A_{d-(i-1)}}n^{d -
(i-1)}\\ &\hskip2in + \mu_{A_{d-i}}n^{d-i} + \mu_{A_{d-(i+1)}}n^{d-(i+1)} + ... +
\mu_{A_1}n + \mu_{A_0}\\ 
&= \mu_{(H + nD_1) \cdot (H + nD_2) \cdot ... \cdot (H +
nD_{d-1})}({\Cal F}). \\
\endalign$$
Therefore, in order for ${\Cal F}$ to have the maximal (averaged) slope with respect to
\linebreak
$(H + nD_1) \cdot (H + nD_2) \cdot ... \cdot (H + nD_{d-1})$ for $n \geq n_i$ it is
necessary that
$$\mu_{A_{d-1}}({\Cal F}) = m_{d-1}, \mu_{A_{d-2}}({\Cal F}) = m_{d-2}, ... ,
\mu_{A_{d-(i-1)}}({\Cal F}) = m_{d-(i-1)}, \mu_{A_{d-i}}({\Cal F}) = m_{d-i}.$$

Finally we conclude that for $n \geq n_d$ in order for a subsheaf ${\Cal F} \subset
{\Cal E}$ to have the maximal (averaged) slope with respect to
\linebreak
$(H + nD_1) \cdot (H + nD_2) \cdot ... \cdot (H + nD_{d-1})$ it is
necessary (and sufficient) that
$$\mu_{A_{d-1}}({\Cal F}) = m_{d-1}, \mu_{A_{d-2}}({\Cal F}) = m_{d-2}, ... ,
\mu_{A_1}({\Cal F}) = m_1, \mu_{A_0}({\Cal F}) = m_0.$$
Since ${\Cal E}_1$ is the maximal subsheaf among such, it is the maximal destabilizing
subsheaf of ${\Cal E}$ with respect to $(H + nD_1) \cdot (H + nD_2) \cdot ... \cdot (H +
nD_{d-1})$ for $n \geq n_d$.

\vskip.1in

This concludes the proof of Proposition 6.4.

\vskip.1in

\proclaim{Proposition 6.5} Let $D_1, D_2, ... , D_{d-1}$ be nef Cartier divisors on
$X$.  Let ${\Cal E}$ be a $(D_1 \cdot D_2 \cdot ... \cdot D_{d-1})$-semistable
torsion free coherent sheaf on $X$.  Suppose that $D_1 \cdot D_2 \cdot ... \cdot
D_{d-1}$ is not numerically trivial and that $\hat{c}_1({\Cal E}) \cdot D_1 \cdot D_2
\cdot ... \cdot D_{d-1} = 0$.  Then
$$\hat{c}_2({\Cal E}) \cdot D_1 \cdot D_2 \cdot ... \cdot D_{d-2} \geq
\frac{1}{2}\hat{c}_1({\Cal E})^2 \cdot D_1 \cdot D_2 \cdot ... \cdot D_{d-2}.$$
In particular, if we assume further that $\hat{c}_1({\Cal E})^2 \cdot D_1 \cdot D_2
\cdot ... \cdot D_{d-2} \geq 0$, then we have
$$\hat{c}_2({\Cal E}) \cdot D_1 \cdot D_2 \cdot ... \cdot D_{d-2} \geq 0.$$
\endproclaim

\demo{Proof}\enddemo By Fujita[3], we take a very ample divisor $H$ on $X$ such that $H
+ nD$ is very ample for any nef Cartier divisor on $X$ and any $n \in {\Bbb N}$.  (The
use of a theorem of Fujita here is purely for the simplicity of presentation.  We could
take a sufficiently high multiple $m(H + nD)$ to make it very ample, and then divide it
back by $m$ in the computation without using the theorem.)  By Proposition 6.4, take
$n_o
\in {\Bbb N}$ such that for
$n \geq n_o$ the Harder-Narasimhan filtration of ${\Cal E}$ with respect to $(H + nD_1)
\cdot (H + nD_2) \cdot ... \cdot (H + nD_{d-1})$
$$0 = {\Cal E}_{0,n} \subset {\Cal E}_{1,n} \subset ... \subset {\Cal E}_{k_n,n} =
{\Cal E}$$
stabilizes to
$$0 = {\Cal E}_0 \subset {\Cal E}_1 \subset ... \subset {\Cal E}_k =
{\Cal E}.$$
Let $S_n = \cap_{i = 1}^{d-2}W_{i,n}$ be the complete intersection of the general
members \linebreak
$W_{i,n} \in |H + nD_i|$.  Then by Mehta-Ramanathan[12], the restriction to $S_n$ of the
Harder-Narasimhan filtration above
$$0 = {\Upsilon}_0 = {\Cal E}_0|_{S_n} \subset {\Upsilon}_1 = {\Cal E}_1|_{S_n}
\subset ... \subset {\Upsilon}_k = {\Cal E}_k|_{S_n} = {\Cal E}|_{S_n} = \Upsilon$$
is the Harder-Narasimhan filtration of ${\Cal E}|_{S_n} = \Upsilon$ on $S_n$.  

Set ${\Cal G}_i =
{\Cal E}_i/{\Cal E}_{i-1}$ and set $g_i = {\Cal G}_i|_{S_n}$, which is a $(H +
nD_{d-1})$-semistable sheaf on $S_n$.

We compute
$$\align
&\hat{c}_2({\Cal E}) \cdot (H + nD_1) \cdot (H + nD_2) \cdot ... \cdot (H + nD_{d-2}) \\
&= \hat{c}_2(\Upsilon) \\
&= \prod_{1 \leq i < j \leq
k}\hat{c}_1(\Upsilon_i/\Upsilon_{i-1})\hat{c}_1(\Upsilon_j/\Upsilon_{j-1}) +
\Sigma_{i=1}^k\hat{c}_2(\Upsilon_i/\Upsilon_{i-1}) \\
&= \prod_{1 \leq i < j \leq
k}\hat{c}_1(g_i)\hat{c}_1(g_j) +
\Sigma_{i=1}^k\hat{c}_2(g_i) \\
&= \frac{1}{2}\hat{c}_1(\Upsilon)^2 +
\Sigma_{i=1}^k\hat{c}_2(g_i) - \frac{1}{2}\Sigma_{i = 1}^k\hat{c}_1(g_i)^2 \\
&\geq \frac{1}{2}\hat{c}_1(\Upsilon)^2 +
\Sigma_{i=1}^k\frac{r_i - 1}{2r_i}\hat{c}_1(g_i)^2 - \frac{1}{2}\Sigma_{i =
1}^k\hat{c}_1(g_i)^2 \\
&= \frac{1}{2}\hat{c}_1({\Cal E})^2 \cdot (H + nD_1) \cdot (H + nD_2) \cdot ... \cdot
(H + nD_{d-2}) - \Sigma_{i = 1}^k\frac{1}{2r_i}\hat{c}_1(g_i)^2 \\
\endalign$$
where $r_i = \roman{rank}\ g_i$ and we used the Bogomolov inequality (cf. Lemma 10.11 in
Koll\'ar et al[9]) to derive the inequality in the second last line.

Moreover, we have by the Hodge index theorem
$$\align
\hat{c}_1(g_i)^2 &\leq \frac{\{\hat{c}_1(g_i) \cdot (H + nD_{d-1}|_{S_n})\}^2}{(H +
nD_{d-1}|_{S_n})^2} \\
&= \frac{\{\hat{c}_1({\Cal G}_i) \cdot (H + nD_1) \cdot (H + nD_2) \cdot ... \cdot (H +
nD_{d-2}) \cdot (H + nD_{d-1})\}^2}{(H + nD_1) \cdot (H + nD_2) \cdot ... \cdot (H +
nD_{d-2})
\cdot (H + nD_{d-1})^2}.\\
\endalign$$

We claim that
$$\hat{c}_1({\Cal
G}_i) \cdot (H + nD_1) \cdot (H + nD_2) \cdot ... \cdot (H + nD_{d-1})$$
is $O(n^{d-2})$, i.e., its absolute value is bounded from above by $cn^{d-2}$, where
$c$ is some positive fixed constant independent of $n$.  This is equivalent to claiming
$$\hat{c}_1({\Cal
G}_i) \cdot D_1 \cdot D_2 \cdot ... \cdot D_{d-1} = 0.$$ 

Postponing the proof of the
claim till the end, we finish the proof of the asserted inequality.  We compute
$$\align
&\hat{c}_2({\Cal E}) \cdot D_1 \cdot D_2 \cdot ... \cdot D_{d-2} \\
&= \lim_{n \rightarrow \infty}\frac{1}{n^{d-2}}\hat{c}_2({\Cal E}) \cdot (H + nD_1)
\cdot (H + nD_2) \cdot ... \cdot (H + nD_{d-2}) \\
&\geq \lim_{n \rightarrow \infty}
\frac{1}{n^{d-2}}[\frac{1}{2}\hat{c}_1({\Cal E})^2
\cdot (H + nD_1)
\cdot (H + nD_2) \cdot ... \cdot (H + nD_{d-2}) \\
&- \Sigma_{i = 1}^k\frac{1}{2r_i}\frac{\{\hat{c}_1({\Cal G}_i) \cdot (H + nD_1) \cdot
(H + nD_2) \cdot ... \cdot (H + nD_{d-2}) \cdot (H + nD_{d-1})\}^2}{(H + nD_1) \cdot (H
+ nD_2) \cdot ... \cdot (H + nD_{d-2})
\cdot (H + nD_{d-1})^2}]\\
&= \frac{1}{2}\hat{c}_1({\Cal E})^2 \cdot D_1 \cdot D_2 \cdot ... \cdot D_{d-2},\\
\endalign$$
obtaining the desired inequality.  

Remark that in concluding
$$\lim_{n \rightarrow \infty}\frac{1}{n^{d-2}}\frac{\{\hat{c}_1({\Cal G}_i) \cdot (H +
nD_1) \cdot (H + nD_2) \cdot ... \cdot (H + nD_{d-2}) \cdot (H + nD_{d-1})\}^2}{(H +
nD_1) \cdot (H + nD_2) \cdot ... \cdot (H + nD_{d-2})
\cdot (H + nD_{d-1})^2} = 0,$$
we used the claim as well as
$$H \cdot D_1 \cdot D_2 \cdot ... \cdot D_{d-1} > 0$$
since $D_1 \cdot D_2 \cdot ... \cdot D_{d-1}$ is not numerically trivial and $H$ is
(very) ample.

\vskip.1in

It remains to prove the claim, by induction on $i$.

\vskip.1in

For $i = 0$, we have ${\Cal G}_1 = {\Cal E}_1$.  Then
$$\align
&\hat{c}_1({\Cal G}_1) \cdot D_1 \cdot D_2 \cdot ... \cdot D_{d-1} \\
&= \hat{c}_1({\Cal
E}_1) \cdot D_1 \cdot D_2 \cdot ... \cdot D_{d-1} \\
&= (\roman{rank}\ {\Cal E}_1) \cdot \mu_{A_{d-1}}({\Cal E}_1) = (\roman{rank}\ {\Cal
E}_1) \cdot m_{d-1}\\
\endalign$$
where $m_{d-1} = \max\{\mu_{A_{d-1}}{\Cal F};{\Cal F} \subset {\Cal E}\}$ and
where $A_{d-1} = D_1 \cdot D_2 \cdot ... \cdot D_{d-1}$. 

(See the proof of Proposition
6.4.)  

Now since
$$\hat{c}_1({\Cal E}) \cdot D_1 \cdot D_2 \cdot ... \cdot D_{d-1} = 0$$
and since ${\Cal E}$ is $A_{d-1} = (D_1 \cdot D_2 \cdot ... \cdot D_{d-1})$-semistable,
we conclude that $m_{d-1} = 0$ and hence that
$$\hat{c}_1({\Cal G}_1) \cdot D_1 \cdot D_2 \cdot ... \cdot D_{d-1} = \hat{c}_1({\Cal
E}_1) \cdot D_1 \cdot D_2 \cdot ... \cdot D_{d-1} = 0.$$
Note that this also implies the quotient ${\Cal E}/{\Cal E}_1$ is $A_{d-1} = (D_1 \cdot
D_2 \cdot ... \cdot D_{d-1})$-semistable with
$$\hat{c}_1({\Cal E}/{\Cal E}_1) \cdot D_1 \cdot D_2 \cdot ... \cdot D_{d-1} = 0.$$

Inductively, we assume that
$$\hat{c}_1({\Cal G}_j) \cdot D_1 \cdot D_2 \cdot ... \cdot D_{d-1} = 0\text{\ for\ }j =
1, ..., i-1$$
and hence
$$\hat{c}_1({\Cal E}_{i-1}) = \Sigma_{j = 1}^{i-1} \hat{c}_1({\Cal G}_j) \cdot D_1
\cdot D_2 \cdot ... \cdot D_{d-1} = 0$$
and that ${\Cal E}/{\Cal E}_{i-1}$ is $A_{d-1} = (D_1 \cdot
D_2 \cdot ... \cdot D_{d-1})$-semistable with
$$\hat{c}_1({\Cal E}/{\Cal E}_{i-1}) \cdot D_1 \cdot D_2 \cdot ... \cdot D_{d-1} = 0.$$
From these assumptions we conclude, as in for $i = 0$ and as in the proof of
Proposition 6.4, that the first piece ${\Cal G}_i = {\Cal E}_i/{\Cal E}_{i-1}$ of the
Harder-Narasimhan filtration of
${\Cal E}/{\Cal E}_{i-1}$ with respect to $(H + nD_1) \cdot (H + nD_2) \cdot ... \cdot
(H + nD_{d-1})$ (for $n \geq n_o$) has the property $\mu_{A_{d-1}}({\Cal G}_i) = 0$ and
hence
$$\hat{c}_1({\Cal G}_i) \cdot D_1 \cdot D_2 \cdot ... \cdot D_{d-1} = 0.$$
Note that this also implies the quotient ${\Cal E}/{\Cal E}_i$ is $A_{d-1} = (D_1 \cdot
D_2 \cdot ... \cdot D_{d-1})$-semistable with
$$\hat{c}_1({\Cal E}/{\Cal E}_i) \cdot D_1 \cdot D_2 \cdot ... \cdot D_{d-1} = 0.$$
Now the inductional argument is complete.

\vskip.1in

This finishes the proof of the claim, and hence that of Proposition 6.5.

\vskip.1in

\proclaim{6.6. Proof of Theorem 6.3}\endproclaim Suppose (2) does not hold.  (That is to
say, suppose either that the tangent sheaf
${\Cal E} = T_X$ is not $(D_1 \cdot D_2 \cdot ... \cdot D_{d-1})$-semistable or that $-
K_X
\cdot D_1
\cdot D_2
\cdot ... \cdot D_{d-2} > 0$.)  Then
$$\hat{c}_1({\Cal E}_{\infty}) \cdot D_1 \cdot D_2 \cdot ... \cdot D_{d-1} > 0,$$ 
where
${\Cal E}_{\infty}$ is the maximal destabilizing subsheaf of ${\Cal E} = T_X$ with
respect to \linebreak
$(D_1 \cdot D_2 \cdot ... \cdot D_{d-1})$ (i.e., the first piece of
the Harder-Narasimhan filtration of ${\Cal E} = T_X$ with respect to $(D_1 \cdot
D_2
\cdot ...
\cdot D_{d-1})$).

By Proposition 6.4, take $n_o \in {\Bbb N}$ such that for
$n \geq n_o$ the Harder-Narasimhan filtration of ${\Cal E}$ with respect to $(H + nD_1)
\cdot (H + nD_2) \cdot ... \cdot (H + nD_{d-1})$
$$0 = {\Cal E}_{0,n} \subset {\Cal E}_{1,n} \subset ... \subset {\Cal E}_{k_n,n} =
{\Cal E}$$
stabilizes to
$$0 = {\Cal E}_0 \subset {\Cal E}_1 \subset ... \subset {\Cal E}_k =
{\Cal E}.$$
Then by the proof of Proposition 6.4, we have
$$\mu_{A_{d-1}}({\Cal E}_1) = m_{d-1} = \mu_{A_{d-1}}({\Cal E}_{\infty}) > 0$$
where $A_{d-1} = (D_1 \cdot D_2 \cdot ... \cdot D_{d-1})$.

Therefore, for $n \gg 0$, we have
$$\hat{c}_1({\Cal E}_1) \cdot (H + nD_1) \cdot (H + nD_2) \cdot ... \cdot (H +
nD_{d-1}) \geq cn^{d-1}$$
where $c$ is some positive constant independent of $n$.

Then, by applying a theorem of Shepherd-Barron (Theorem 9.0.2 in Koll\'ar et al[9])
after Miyaoka-Mori[16], we have a covering family of rational curves $C$ such that
$$(H + nD_d) \cdot C \leq \frac{2d(H + nD_1) \cdot (H + nD_2) \cdot ... \cdot (H +
nD_d)}{\hat{c}_1({\Cal E}_1) \cdot (H + nD_1) \cdot (H + nD_2) \cdot ... \cdot (H +
nD_{d-1})}.$$ 
Remark that
$$\align
&\lim_{n \rightarrow \infty}\frac{2d(H + nD_1) \cdot (H + nD_2) \cdot ... \cdot (H +
nD_d)}{\hat{c}_1({\Cal E}_1) \cdot (H + nD_1) \cdot (H + nD_2) \cdot ... \cdot (H +
nD_{d-1})} \\
&\leq \lim_{n \rightarrow \infty}\frac{2d(H + nD_1) \cdot (H + nD_2) \cdot ... \cdot (H
+ nD_d)}{cn^{d-1}} < \alpha < \infty,\\
\endalign$$
since
$$D_1 \cdot D_2 \cdot ... \cdot D_d = 0.$$
Thus for $n \gg 0$ we have $(H + nD_d) \cdot C < \alpha$, and in particular if we take
$n \geq \alpha$ then
$$D_d \cdot C = 0$$
since both $H \cdot C$ and $D_d \cdot C$ are nonnegative integers.

This finishes the proof of the first part of Theorem 6.3.

\vskip.1in

The second part, i.e., ``Moreover" part of Theorem 6.3, follows directly from
Proposition 6.5.

\vskip.1in

This finishes the proof of Theorem 6.3.

\vskip.1in

We mention two direct corollaries of Theorem 6.3, which are due to S. Keel and J.
McKernan.

\proclaim{6.7 Corollary} Suppose $X$ is a normal projective variety of dimension
$d$ with only canonical singularities.  Suppose
$- K_X$ is nef of numerical dimension
$m < d$.  Assume that $D_1, D_2, ... , D_{d-m-1}$ are nef Cartier
divisors such that \linebreak
$(- K_X)^m \cdot D_1 \cdot D_2 \cdot ... \cdot D_{d-m-1}$
is not numerically trivial.  

Then $\hat{c}_2 \cdot (- K_X)^{m-1} \cdot D_1 \cdot D_2
\cdot ...
\cdot D_{d-m-1}$ is nonnegative.
\endproclaim

\demo{Proof}\enddemo Apply Theorem 6.3 with $D_{d-m} = D_{d-m+1} = ... = D_d = -
K_X$.  Remark that there is no covering family of rational curves $C$ with $- K_X \cdot
C = 0$, and hence case (1) never occurs.

\proclaim{6.8 Corollary} Let $X$ be a normal projective 3-fold with only canonical
singularities.  Suppose that $- K_X$ is nef but not big, and that
$\chi({\Cal O}_X) = 0$.  Then $X$ is Gorenstein.
\endproclaim

\demo{Proof}\enddemo Recall (See (10.3) of Reid [19].) that the Riemann-Roch formula
for a threefold with canonical singularities implies
$$\chi({\Cal O}_X) = \frac{1}{24}\hat{c}_1\hat{c}_2 + x$$
where $x$ is a nonnegative term which is zero if and only if $X$ is Gorenstein.  Thus
it is enough to show $\hat{c}_1\hat{c}_2 = \hat{c}_2 \cdot (- K_X) \geq 0$ under the
assumption.  When $\nu(- K_X) = 0$, we have $\hat{c}_2 \cdot (- K_X) = 0$.  When $\nu(-
K_X) = 1$, take an ample divisor $H$ with $(- K_X) \cdot H$ being not numerically
trivial.  Then apply Theorem 6.3, setting $D_1 = - K_X, D_2 = H$.  When $\nu(- K_X) =
2$, apply Theorem 6.3, setting $D_1 = D_2 = - K_X$.  Remark again that there is no
covering family of rational curves $C$ with $- K_X \cdot C = 0$, and
hence case (1) never occurs.

\vskip.2in

$$\bold{References}$$

[1] F.A. Bogomolov, \it Holomorphic tensors and vector bundles on projective 

\hskip.1in varieties,
\rm Math. U.S.S.R. Izvestja $\bold{13}$ (1979), 499-555

[2] D. Eddidin, B. Hassett, A. Kresch, and A. Vistoli, \it Brauer groups and quotient

\hskip.1in stacks, \rm Amer. J. Math. $\bold{123}$ (2001), 761-777

[3] T. Fujita, \it Semipositive line bundles, \rm J. Fac. Sci. Tokyo Univ. $\bold{30}$
(1983), 

\hskip.1in 353-378

[4] Y. Kawamata, \it Pluricanonical systems on minimal algebraic varieties, \rm Invent.

\hskip.1in Math. $\bold{79}$ (1985), 567-588

[5] Y. Kawamata, \it Abundance theorem for minimal threefolds, \rm Invent. Math.
$\bold{108}$ 

\hskip.1in (1992), 229-246

[6] Y. Kawamata, \it Log canonical models of threefolds, \rm Intern. J. Math. $\bold{3}$
(1992), 

\hskip.1in 351-357

[7] Y. Kawamata, K. Matsuda, and K. Matsuki, \it Introduction to the minimal 

\hskip.1in model problem, \rm Adv. Stud. Pure Math. $\bold{10}$, Alg. Geom.,
Sendai, T. Oda ed. 

\hskip.1in (1987), 283-360

[8] S. Keel, K. Matsuki, and J. McKernan, \it Log Abundance theorem for threefolds,

\hskip.1in \rm Duke Math. J. $\bold{75}$ (1994), 99-119 

[9] Koll\'ar et al, \it Flips and Abundance for Algebraic Threefolds, \rm Ast\'erisque
$\bold{211}$ 

\hskip.1in (1993)

[10] K. Matsuki and R. Wentworth, \it Mumford-Thaddeus principle on the moduli 

\hskip.1in space of
vector bundles on an algebraic surface, \rm Internat. J. Math. $\bold{8}$ (1977), 

\hskip.1in 97-148

[11] T. Matsusaka, \it Theory of ${\Bbb Q}$-varieties, \rm Publication of the
Mathematical 

\hskip.1in Society of Japan $\bold{8}$ (1964) 

[12] V.B. Mehta and A. Ramanathan, \it Semistable sheaves on projective varieties 

\hskip.1in and
their restriction to curves, \rm Math. Ann. $\bold{258}$ (1982), 213-224

[13] Y. Miyaoka, \it The Chern classes and Kodaira dimension of a minimal variety, 

\hskip.1in \rm Adv. Stud. Pure Math. $\bold{10}$, Alg. Geom., Sendai, T. Oda ed.
(1987), 449-476 

[14] Y. Miyaoka, \it On the Kodaira dimension of minimal threefolds, \rm Math. Ann.

\hskip.1in $\bold{281}$ (1988), 325-332

[15] Y. Miyaoka, \it Abundance Conjecture for 3-folds: case $\nu = 1$, \rm Compositio

\hskip.1in Math. $\bold{68}$ (1988), 203-220

[16] Y. Miyaoka and S. Mori, \it A numerical criterion for uniruledness, \rm Ann. of

\hskip.1in Math. $\bold{124}$ (1986), 65-69

[17] S. Mori, \rm Projective manifolds with ample tangent bundles, \rm Ann. of Math.

\hskip.1in $\bold{110}$ (1986), 593-606 

[18] D. Mumford, \it Towards an enumerative geometry of the moduli space of 

\hskip.1in curves,
\rm Arithmetic and Geometry, Papers dedicated to I.R. Shafarevich on the 

\hskip.1in occasion of his
60th birthday, Volume II, Geometry, M. Artin and J. Tate eds., 

\hskip.1in Birkh\"auser (1983),
271-328

[19] M. Reid, \it Young person's guide to canonical singularitiers, \rm in Algebraic

\hskip.1in Geometry, Bowdoin, 1985, Proc. Symp. Pure Math. $\bold{46}$, Amer. Math. Soc.

\hskip.1in (1987), 345-416

[20] B. Totaro, \it The resolution property for schemes and stacks, \rm math.AG/0207210

[21] A. Vistoli, \it Intersection theory on algebraic stacks and on their moduli spaces,

\hskip.1in \rm Invent. Math. $\bold{97}$ (1989), 613-670

[22] Q. Xie, \it private communication, \rm (2002)

\enddocument